\documentclass[runningheads]{lncse}

\usepackage{amsmath}
\usepackage{amsfonts}
\usepackage{epsfig}
\usepackage{graphics} 
\usepackage{subfigure}

\renewcommand{\phi}{\varphi}

\newcommand{\trans}{^{\sf T}}

\def\~{\widetilde}
\def\^{\widehat}

\newcommand{\ii}{\text{i}\hspace{1pt}}
\newcommand{\dd}{\hspace{1pt}{\rm d}\hspace{0.5pt}}

\newcommand{\RR}{\mathbb{R}}

\newcommand{\CC}{\mathbb{C}}

\newcommand{\NN}{\mathbb{N}}

\newcommand{\veps}{\varepsilon}

\begin{document}

\title{On a generalization of the Bessel function Neumann expansion}

\titlerunning{Generalization of Bessel-Neumann expansion}

\author{A. Koskela\inst{1} \and E. Jarlebring\inst{1}}

\authorrunning{A. Koskela, E. Jarlebring}   

\institute{
KTH Royal Institute of Technology {\em [akoskela, eliasj]@kth.se} 
}

\maketitle

\begin{abstract}
  The Bessel-Neumann expansion (of integer order) of a function $g:\CC\rightarrow\CC$
  corresponds to representing $g$ as
  a linear combination of basis functions $\phi_0,\phi_1,\ldots$, i.e.,
  $g(z)=\sum_{\ell = 0}^\infty w_\ell \phi_\ell(s)$, where $\phi_i(z)=J_i(z)$, $i=0,\ldots$,
  are the Bessel functions.
  In this work, we study an expansion for a more general class of basis functions.
  More precisely, we assume that the basis functions
  satisfy an infinite dimensional linear ordinary differential equation
  associated with a Hessenberg matrix, motivated by the fact that these
  basis functions occur in certain iterative methods.  A procedure 
  to compute the basis functions as well as the coefficients is proposed.
  Theoretical properties of the expansion are studied. 
  We illustrate that non-standard basis functions can give faster convergence than the Bessel functions.
\end{abstract}

\section{Introduction}
Let $g:\CC\rightarrow\CC$ be a function analytic in a
neighborhood of the origin containing $s \in \mathbb{R}$.
We consider expansions of the form 
\begin{equation} \label{eq:expansion}
g(s) = \sum\limits_{\ell = 0}^\infty w_\ell \phi_\ell(s),
\end{equation}
where the coefficients $w_0,w_1,\ldots\in\CC$ and the basis
functions $\phi_0,\phi_1,\ldots$ satisfy
\begin{equation}  \label{eq:inf_system}
\frac{\dd}{\dd t} \begin{bmatrix} \phi_0(t) \\ \phi_1(t) \\ \vdots \end{bmatrix} = H_\infty 
\begin{bmatrix} \phi_0(t) \\ \phi_1(t) \\ \vdots \end{bmatrix}, \quad \quad \begin{bmatrix} \phi_0(0) \\ \phi_1(0) \\ \vdots \end{bmatrix}=
e_1.
\end{equation}
The (infinite) matrix $H_\infty$ is an infinite upper Hessenberg matrix with non-zero subdiagonal entries. Due
to the Picard-Lindel\"of theorem, the
differential equation \eqref{eq:inf_system} defines
a unique sequence of functions $\phi_0,\phi_1,\ldots$
under the condition that 
\begin{equation}  \label{eq:boundedness}
\|H_\infty\|_2\le C<\infty.
\end{equation}
This implies that also $\|H_n\|_2\le C$, where $H_n\in\RR^{n\times n}$ is the leading submatrix of $H_\infty$.


For certain choices of $H_\infty$
\eqref{eq:inf_system} defines well-known basis functions. We will show
in Section~\ref{sect:wellknown} that this class of basis functions includes
monomials, Bessel functions and modified Bessel functions (of the first kind) and
the expansion corresponds to well-known expansions (e.g. Taylor, Bessel function Neumann expansion).
Although our work is of general character, it
originated from the need of such expansions in iterative methods.
More precisely, the expansion occurs in the infinite Arnoldi method
both for linear ODEs \cite{Koskela:2016:ODEpoly,Koskela:2015:SEMILIN} and
similarly for nonlinear
eigenvalue problems \cite{Jarlebring:2010:DELAYARNOLDI},
where the fact that $H_\infty$ is a Hessenberg matrix leads to the
property that the $m$ does not have to be determined a priori.


In this work we study this expansion, both from a computational and
theoretical perspective.
\begin{itemize}
  \item  We derive procedures to compute
    the coefficients $w_0,w_1,\ldots,w_{n-1}$, based on the derivatives of $g$ and $H_n$.
   \item  We derive a procedure to compute the basis functions 
based on a	 truncation of $H_\infty$  using the matrix exponential \cite{Moler:2003:NINETEEN}, and we provide
a convergence analysis of the corresponding problem.
\item We provide an error analysis for the truncation of the series, i.e., for the error of the the approximation
\begin{equation} \label{eq:truncexpansion}
g(s) \approx \sum\limits_{\ell = 0}^{n-1} w_\ell \phi_\ell(s).
\end{equation}
\end{itemize}


The results are illustrated by comparing different choices of the basis functions, i.e.,
different choices of the matrix $H_\infty$, for a given a function $g(s)$.





\section{Well-known basis functions}\label{sect:wellknown}
We first illustrate that the class of expansions \eqref{eq:expansion} 
includes several well-known expansions. 
The scaled monomials form the simplest example
of such a sequence of basis functions.
If we define $\phi_\ell(t):=t^\ell/\ell!$, $i=0,\ldots$, 
then \eqref{eq:inf_system} is satisfied with $H_\infty$ given by a transposed Jordan matrix
\begin{equation} \label{eq:H_jordan}
H_\infty = \begin{bmatrix}
   0 & & & \\
   1 &0 & & \\
    &1 &0 & & \\
   & &\ddots & \ddots &
 \end{bmatrix}.
\end{equation}
 In this case, the
expansion \eqref{eq:expansion} corresponds to a Taylor expansion and
the coefficients are given by $w_\ell=g^{(\ell)}(0)$, $\ell=0,\ldots$.

It turns out that also Bessel functions and modified Bessel functions
of the first kind satisfy \eqref{eq:inf_system}. It can be seen as follows.
The Bessel functions of the first kind are defined by (see e.g.~\cite{Stegun:1964:HANDBOOK})
$J_\ell(t):=\frac{1}{\pi}\int_{0}^\pi\cos(\ell\tau -t\sin(\tau))\,d\tau$,  
for $\ell\in\NN$, and they satisfy 
\begin{equation}\label{eq:Bessel_properties}
\begin{aligned}
J_\ell'(t)&=\frac12 (J_{\ell-1}(t)-J_{\ell+1}(t)).\\
J_{-\ell}(t)& = (-1)^\ell J_\ell(t),\;\;\ell>0\\
J_\ell(0)&=
\begin{cases}
   1 & \textrm{ if } \ell=0\\
   0 & \textrm{ otherwise}\\
\end{cases}.
\end{aligned}
\end{equation}
Let $\bar J_n(t) = \begin{bmatrix} J_0(t) & J_1(t) & \dots & J_{n-1}(t) \end{bmatrix}^\mathrm{T}   \in\RR^{n}$,
i.e., a vector of Bessel functions
with non-negative index. Moreover, let
\begin{equation} \label{eq:H_bessel}
 H_n = 
\begin{bmatrix} 0 & -1 & &  & \\
	      \tfrac{1}{2} & 0 & -\tfrac{1}{2} &  &  \\
	         & \ddots & \ddots & \ddots &  \\
	         & &\ddots & 0 & -\tfrac{1}{2} &  \\
	         & & &\tfrac{1}{2} &0 & \end{bmatrix}
\in\CC^{n\times n}.
	         \end{equation}
From the  relations \eqref{eq:Bessel_properties}, 
we easily verify that the Bessel functions of the first kind are 
solutions to the infinite-dimensional ODE of the form \eqref{eq:inf_system},
with $H_\infty$ defined by $\eqref{eq:H_bessel}$. More precisely, 
$$
\bar J_\infty'(t) = H_\infty \bar J_\infty(t), \quad \bar J_\infty(0) = e_1.
$$
%
With similar reasoning we can establish an 
ODE \eqref{eq:inf_system} also for  the modified Bessel functions of the first kind, which are defined by
\begin{equation} \label{eq:modified_Bessels}
I_\ell(t) := (- \ii)^\ell J_\ell(\ii t).
\end{equation}
The definition \eqref{eq:modified_Bessels} and the properties \eqref{eq:Bessel_properties} lead to the infinite-dimensional ODE 
$$
\bar I_\infty'(t) = H_\infty \bar I_\infty(t), \quad I(0) = e_1,
$$
where $\bar I_n(t) = \begin{bmatrix} I_0(t)  & I_1(t) & \dots & I_{n-1}(t) \end{bmatrix}\trans$ and 
\begin{equation} \label{eq:H_bessel_modified}
 H_n = 
\begin{bmatrix} 0 & 1 & &  & \\
	      \tfrac{1}{2} & 0 & \tfrac{1}{2} &  &  \\
	         & \ddots & \ddots & \ddots &  \\
	         & &\ddots & 0 & \tfrac{1}{2} &  \\
	         & & &\tfrac{1}{2} &0 & \end{bmatrix}
\in\CC^{n\times n}.
\end{equation}
Therefore the Bessel functions and the modified Bessel functions  of the first kind satisfy
\eqref{eq:inf_system}. The Euclidean norm bound
of the operator ($C$) can be derived for these choices as proven with
recurrence relations of Chebyshev polynomials in \cite[Lemma~2]{Koskela:2015:SEMILIN}.
\begin{lemma}[Basis functions. Lemma~2 in \cite{Koskela:2015:SEMILIN}]
The conditions for the basis functions 
in \eqref{eq:inf_system} are satisfied with $C=2$ for,
\begin{itemize}
  \item[(a)] scaled monomials, i.e., $\phi_i(t)=t!/i!$, with $H_\infty$ defined by \eqref{eq:H_jordan};
  \item[(b)] Bessel functions, i.e., $\phi_i(t)=J_i(t)$, with $H_\infty$ defined by \eqref{eq:H_bessel}; and 
  \item[(c)] modified Bessel functions, i.e., $\phi_i(t)=I_i(t)$, with $H_\infty$ defined by \eqref{eq:H_bessel_modified}.
\end{itemize}
\end{lemma}
%

\section{Computation of coefficients $\mathbf{w_i}$}\label{sect:coeffs}
Assume that an expansion of the form \eqref{eq:expansion} exists and let  
$$
W_n = \begin{bmatrix} w_0 & w_1 & \ldots & w_{n-1} \end{bmatrix}.
$$
By considering the $\ell$th derivative of $g(t)$  and
using the properties of basis functions \eqref{eq:inf_system} we have that
\begin{equation} \label{eq:property}
g^{(\ell)}(0) = W_\infty H_\infty^\ell e_1=
W_nH_n^\ell e_1 \quad \textrm{for all} \quad \ell< n.
\end{equation}
In the last equality we used the fact that $H_\infty$ is a Hessenberg matrix, and
that all elements of $H_\infty^\ell e_1$ except the first $\ell+1$ elements will be zero. 
The non-zero elements will also be equal to $H_n^\ell e_1$.
We now define the upper-triangular Krylov matrix 
\begin{equation*} 
K_n(H_n,e_1) = \begin{bmatrix} e_1 & H_n e_1 & \ldots & H_n^{n-1} e_1 \end{bmatrix},
\end{equation*}
and the matrix $G_n$ as
\begin{equation*} 
G_n = \begin{bmatrix} g(0) & g'(0) & \ldots & g^{(n-1)}(0) \end{bmatrix}.
\end{equation*}
From these definitions and the property \eqref{eq:property} it follows that
\begin{equation} \label{eq:Bessel_Krylov_relation}
W_n= G_n  K_n(H_n,e_1)^{-1}    \quad \textrm{for all} \quad n\geq 1.
\end{equation}
Notice that $K_n(H_n,e_1)$ is invertible if and only if the subdiagonal elements of $H_n$ are non-zero.
In a generic situation, the relation \eqref{eq:Bessel_Krylov_relation} 
can be directly used to  compute the coefficients $w_\ell$, $\ell\in\NN$, given the derivatives of $g$.

The above approach reduces to standard procedures of evaluating coefficients in
Bessel-Neumann series (for integer order) when $\phi_i=J_i$ or $\phi_i=I_i$.
The formalization is omitted due to space limitation, but
can be found in the technical report \cite[Lemma~3]{Koskela:2015:SEMILIN}. 
The relation with Chebyshev polynomials (in particular \cite[pp.\;775]{Stegun:1964:HANDBOOK}
and \cite{Watson:1995:TREATISE})
is pointed out in \cite[Remark~4]{Koskela:2015:SEMILIN}.
See also other works on the Bessel-Neumann series coefficients in \cite{Dragana:2011:NEUMANN},
not necessarily of integer order.

\section{Computation of basis functions}\label{sect:basis}
In order to use the expansion in practice, we
need the possibility to compute the basis functions. 
We propose and study a natural approach of approximating
the basis functions with a vector of functions generated by the truncated Hessenberg matrix $H_n$, i.e.,
\begin{equation}
\bar\varphi_n(t):=
\begin{bmatrix}
  \phi_0(t)\\
  \vdots\\
  \phi_{n-1}(t)
 \end{bmatrix}
\approx \exp(tH_n)e_1.
\end{equation}
The term $\exp(tH_n)$ can be efficiently and accurately
computed using well-known methods for the matrix exponential \cite{Moler:2003:NINETEEN}.
Let $\bar\veps_n$ denote the error in the
basis functions, i.e.,
the difference between the basis functions and the functions generated by the truncated Hessenberg matrix,
\begin{equation}\label{eq:vepsdef}
\begin{aligned}
\bar\veps_n(t) := &\bar\phi_n(t)-\exp(tH_n)e_1,\;\;\; \\
= & [I_n, 0]\exp( t H_\infty )e_1-\exp(tH_n)e_1,
\end{aligned}
\end{equation}
where $[I_n, 0]$ gives the first $n$ rows of an (infinite) identity matrix.
Since $H_\infty$ and $H_n$ are Hessenberg matrices, for $i=0,\ldots,n-1$ we have
\[
[I_n, 0] H_\infty^ie_1=H_n^ie_1
\]
such that $n$ first terms in the Taylor expansions of
the matrix exponentials in \eqref{eq:vepsdef} cancel. We can
therefore explicitly bound the basis function error as
\begin{equation*}
\begin{aligned}
  \|\bar\veps_n(t)\|&=  \|[I_n, 0] R_n(tH_\infty)e_1-R_n(tH_n)e_1\| \\
  &\le R_n(t\|H_\infty\|)+R_n(t(\|H_n\|))\\
  &\le  2R_n(t\|H\|)\le \frac{2(tC)^n}{n!}e^{tC}
\end{aligned}
\end{equation*}
where $R_n(z) = \sum_{\ell=n}^\infty \tfrac{z^\ell}{\ell !}$, i.e., the remainder of the truncated Taylor expansion.
This shows that the basis function approximations converge superlinearly as $n \rightarrow\infty$, for a fixed value of $t$. 


Sharper bounds can be derived for special cases.
Sharper results for Bessel functions can be found
in~\cite[Section~4.1.1-4.1.2]{Koskela:2015:SEMILIN}.

\section{Truncation of the expansion}

We now provide sufficient conditions for the convergence of the expansion by giving bounds for
the truncation error
\[
e_n :=g(s)-W_n\bar\varphi_n(s).
\]

\subsection{Truncation bound with matrix exponential remainder term}
The truncation error can now be expressed as infinite matrices as follows:
\begin{equation}\label{eq:en_expm_expression}
e_n=\sum_{p=n}^\infty w_p\varphi_p(s)= W_\infty
\begin{bmatrix}
    0_n&0\\
    0 & I
\end{bmatrix} \exp(sH_\infty)e_1.
\end{equation}
The fact that $H_\infty$ is a Hessenberg matrix implies that $e_j^TH_\infty^ie_1=0$,
when $j>i+1$. Therefore
\[
e_n=\sum_{i=0}^\infty W_\infty
\begin{bmatrix}
    0_n&0\\
    0 & I
\end{bmatrix}\frac{(sH_\infty)^i}{i!}e_1=
W_\infty \begin{bmatrix}
    0_n&0\\
    0 & I
\end{bmatrix}
R_n(sH_\infty)e_1
\]
A simple sufficient condition for the convergence of the series
is obtained by assuming that the sequence $|w_0|,|w_1|,\ldots$ is bounded, since then
\begin{equation*}
\begin{aligned}
|e_n|\le& \|W_\infty\|_2\|R_n(sH_\infty)e_1\|_2\le \|W_{\infty}\|_1R_n(|s|\|H_\infty\|_2)\\
\le& \|W_\infty\|_1R_n(|s|C)\rightarrow 0, \quad \textrm{as} \quad n \rightarrow \infty.
\end{aligned}
\end{equation*}
The above reasoning can be generalized by scaling with a bounding sequence and leads to the following theorem.
The proof follows analogous to the derivation of \eqref{eq:en_expm_expression}.
\begin{theorem}\label{thm:simplebound}
  Suppose the $w_0,w_1,\ldots$ is bounded in modulus by $d_0,d_1,\ldots\in\RR_+$,
  \[
   |w_i|\le d_i,\;\;i=0,1,\ldots.
  \] 
  Let $D= \mathrm{diag}(d_0,d_1,\ldots)$.  Then,
  \begin{equation*} 
  \begin{aligned}
    |e_n| \le & \|D R_n(sH_\infty)e_1\|_2\\
     = &d_0\|R_n(sDH_\infty D^{-1})e_1\|_2.
  \end{aligned}
  \end{equation*}
\end{theorem}
As a consequence of the theorem, if we have 
\[
  \|DH_\infty D^{-1}\|_2\le \infty,
\]
we may use the bound and obtain
  \begin{equation*} 
  \begin{aligned}
  e_n\le & |d_0|\max_i \frac{|w_i|}{|d_i|}R_n(|s|\|D H_\infty D^{-1}\|_2) \\
  \le & |d_0|\max_i \frac{|w_i|}{|d_i|}  e^{|s|\|DH_\infty D^{-1}\|_2}\frac{|s|^n\|DH_\infty D^{-1}\|_2^n}{n!}.
  \end{aligned}
  \end{equation*}  

  \subsection{Bounds using decay of elements in the matrix exponential}
  Equation \eqref{eq:en_expm_expression} can be reformulated as
\begin{equation}\label{eq:en_expm_sum}
\begin{aligned}
  |e_n|  \le  \sum_{j=n+1}^\infty |w_j||e_j^T\exp(sH_\infty)e_1|
\end{aligned}
\end{equation}
or with scaling as in Theorem~\ref{thm:simplebound},
\begin{equation}\label{eq:en_expm_sum_D}
  |e_n|\le \sum_{j=n+1}^\infty d_0|e_j^T\exp(sDH_\infty D^{-1})e_1|.
\end{equation}
Note that $DH_\infty D^{-1}$ is banded if $H_\infty$ is banded. 
The elements of matrix functions of banded matrices are bounded in a ridge shaped way. 
There are bounds on the absolute value of elements which show that elements far from the diagonal
have to be small. These can be used here to analyze 
the bounds \eqref{eq:en_expm_sum} and \eqref{eq:en_expm_sum_D} as they consist of
elements $(j,1)$ of the matrix exponential.

%
%

As an example we mention a result by Iserles (\cite[theorem 2.2]{Iserles:2000:EXPMBAND})
which is applicable for a tridiagonal $H_\infty$ and that of Benzi and Razouk (\cite[theorem 3.5]{Benzi:2007:decay_bounds})
which is also used in~\cite[Lemma 9, Thm.\;11]{Koskela:2015:SEMILIN}.



\section{Illustrating example}
The following example illustrates how a generalized Bessel-Neumann expansion
of the form \eqref{eq:expansion} can lead to a better approximation of the function when the matrix $H_\infty$ is chosen
appropriately. We consider the function 
\begin{equation}\label{eq:g_example}
  g(s)=e^{\alpha s}\left(\sin\left(s/3\right)+\cos(s)\right)
  \end{equation}
  with $\alpha=1/2$. We use the Hessenberg matrices
  corresponding to scaled monomials \eqref{eq:H_jordan},
  Bessel functions  \eqref{eq:H_bessel},
  modified Bessel functions \eqref{eq:H_bessel_modified},
  and an artifically constructed Hessenberg matrix
  given by $H+\alpha I$ where $H$ is given by 
  \eqref{eq:H_bessel}. The infinite matrix $H+\alpha I$ is bounded in Euclidean norm
  since
  $\|H+\alpha I\|\le \|H\|+|\alpha|$. 
  The coefficients and basis functions
  are computed as described in Section~\ref{sect:coeffs}
  and Section~\ref{sect:basis}.  In Figures 1 and 2
  we clearly see that the generalized Bessel-Neumann expansion
  converges faster for this (specific) example than the other alternatives.
  
  \begin{figure}[h!]
  \begin{center}
    \includegraphics[scale=0.42]{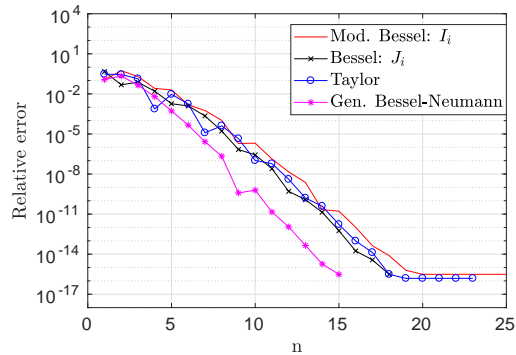}
      \end{center}
    \caption{Relative truncation error as a function of truncation parameter,
      when $s=1$ for the function \eqref{eq:g_example}.  
    }\label{fig:exmp0}
\end{figure}

 \begin{figure}[h!] 
  \begin{center}
    \includegraphics[scale=0.42]{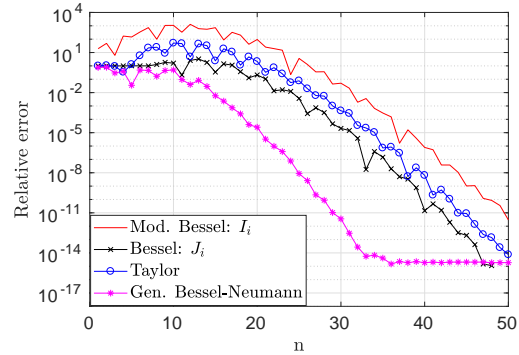}
      \end{center}
    \caption{Relative truncation error as a function of truncation parameter,
      when $s=10$ for the function \eqref{eq:g_example}. 
    } \label{fig:exmp1}
\end{figure}

\section{Concluding remarks}

We have here proposed the use of 
generalized Bessel-Neumann expansions, and shown both
computational and theoretical properties for this approach. The main
purpose of this work is to show that it is possible use these expansions and that the convergence
speed depends on the choice of the matrix $H_\infty$ which determines both the basis functions and 
the coefficients of the expansion. We have illustrated numerically
that different choices of $H_\infty$ lead to different performance, and it is not at all
obvious which basis functions are the best, and in which sense they are good for a particular function.

%
%
%
%
%
%



%
%

%
%
%
\ifx\undefined\bysame
\newcommand{\bysame}{\leavevmode\hbox to3em{\hrulefill}\,}
\fi

\end{document}